\documentclass[11pt,a4paper]{article}

\usepackage{color}
\usepackage{graphics,epsfig}
\usepackage{amsfonts,amssymb,amsmath,amsthm}

\newtheorem{lem}{Lemma}[section]
\newtheorem{cor}[lem]{Corollary}
\newtheorem{thm}[lem]{Theorem}
\newtheorem{prop}[lem]{Proposition}
\theoremstyle{definition}

\newtheorem{defi}[lem]{Definition}
\theoremstyle{remark}

\numberwithin{equation}{section}

\newcommand{\nonlocal}{\gamma}
\newcommand{\f}{\tilde f _\delta}
\newcommand{\de}{\delta}
\newcommand{\ep}{\varepsilon}
\newcommand{\ue}{u^\ep}
\newcommand{\om}{\Omega}
\newcommand{\ombar}{\overline{\Omega}}
\newcommand{\edeux}{\displaystyle{\frac{1}{\ep^2}}}
\newcommand{\Pe}{(P^\ep)}
\newcommand{\R}{\mathbb{R}}
\newcommand{\EB}{e^{-\beta t/\ep ^2}}
\newcommand{\vsp}{\vspace{8pt}}
\newcommand{\n}{\nabla}
\newcommand{\di}{\displaystyle}
\newcommand{\hm}{\alpha _- (\de)}
\newcommand{\hp}{\alpha _+ (\de)}
\newcommand{\h}{a(\de)}
\newcommand{\mm}{\mu(\de)}
\newcommand{\emutt}{e^{\mu (\ep \mathcal G) t/\ep ^2}}
\newcommand{\cste}{C_\star}
\newcommand{\autrecste}{M_ \star}
\newcommand{\epaisseur}{\mathcal C}

\title{Generation, motion and thickness of transition layers for a nonlocal Allen-Cahn equation}
\author{ }
\date{}

\begin{document}
\maketitle \vspace{-15 mm}

\begin{center}

{\large \bf Matthieu Alfaro}\\[1ex]
D\'epartement de Math\'ematiques, CC051, Universit\'e Montpellier II,\\
Place Eug\`ene Bataillon, 34095 Montpellier Cedex 5, France.\\
email: malfaro@math.univ-montp2.fr \\
fax: +33 (0)4 67 14 93 16\\[2ex]
\end{center}

\vspace{15pt}

\begin{abstract}
We investigate the behavior, as $\ep \to 0$, of the nonlocal
Allen-Cahn equation $u_t=\Delta u+\edeux f(u,\ep \int _ \Omega
u)$, where $f(u,0)$ is of the bistable type. Given a rather
general initial data $u_0$ that is independent of $\ep$, we
perform a rigorous analysis of both the generation and the motion
of interface, and obtain a new estimate for its thickness. More
precisely we show that the solution develops a steep transition
layer within the time scale of order $\ep^2|\ln\ep|$, and that the
layer obeys the law of motion that coincides with the limit
problem within an error margin of order $\ep$.\\

\end{abstract}

\section{Introduction}\label{s:intro}
This paper is concerned with the singular limit, as $\ep \to 0$,
of the nonlocal Allen-Cahn equation
\[
 \Pe \quad\begin{cases}
 u_t=\Delta u+\edeux f\left(u,\ep \int _ \Omega u\right) &\text{in }\om \times (0,\infty)\\
 \displaystyle{\frac{\partial u}{\partial\nu}} = 0 &\text{on }\partial \om \times (0,\infty)\vspace{3pt}\\
 u(x,0)=u_0(x) &\text{in }\om,
 \end{cases}
\]
where $\om$ is a smooth bounded domain in $\R ^N$ ($N\geq 2$) and
$\nu$ the Euclidian unit normal vector exterior to $\partial \om$.
We assume that the nonlinearity $f(u,v)$ is smooth and that
$\tilde f(u):=f(u,0)$ is given by $\tilde f(u):=-W'(u)$, where
$W(u)$ is a double-well potential with equal well-depth, taking
its global minimum value at $u=\pm 1$. More precisely we assume
that $\tilde f$ has exactly three zeros $-1<a<1$ such that
\begin{equation}\label{der-f}
\tilde f '(\pm 1)<0, \quad \tilde f'(a)>0\quad\ \hbox{(bistable
nonlinearity)},
\end{equation}
and that
\begin{equation}\label{int-f}
\int _ {-1} ^ {+1} \tilde f(u)\,du=0.
\end{equation}
The condition \eqref{der-f} implies that the potential $W(u)$
attains its local minima at $u=\pm 1$, and \eqref{int-f} implies
that $W(-1)=W(+1)$. In other words, the two stable zeros of
$\tilde f$ have ``balanced" stability.

Concerning the initial data $u_0$, we assume its smoothness and
choose $C_0\geq 1$ such that
\begin{equation}\label{int1}
\|u_0\|_{C^0(\ombar)}+\|\n u_0\|_{C^0(\ombar)}+\|
D^2u_0\|_{C^0(\ombar)}\leq C_0.
\end{equation}
Furthermore we define the \lq \lq initial interface" $\Gamma _0$
by
\[\Gamma _0:=\{x\in\om| \; u_0(x)=a \},\]
and suppose that $\Gamma _0$ is a smooth closed hypersurface
without boundary, such that, $n$ being the Euclidian unit normal
vector exterior to $\Gamma _0$,
\begin{equation}\label{dalltint-aniso}
\Gamma _0 \subset\subset \Omega \quad \mbox { and } \quad \n
u_0(x) \neq 0\quad\text{if $x\in\Gamma _0,$}
\end{equation}
\begin{equation}\label{initial-data-aniso}
u_0>a \quad \text { in } \quad \om ^+ _0,\quad u_0<a \quad \text {
in } \quad \om ^- _0 ,
\end{equation}
where $\om ^- _0$ denotes the region enclosed by $\Gamma _0$ and
$\om ^+ _0$ the region enclosed between $\partial \om$ and $\Gamma
_0$.

\vskip 8 pt Before going into more details, let us recall known
facts concerning the \lq\lq usual" Allen-Cahn equation, namely
$$
u_t=\Delta u+\edeux \tilde f(u).
$$
The singular limit was first studied by Allen and Cahn \cite{AC}
and by Kawasaki and Ohta \cite{KO}. By using formal asymptotic
arguments, they show that the limit problem, as $\ep \to 0$,  is a
free boundary problem: the motion of the limit interface is ruled
by its mean curvature. More precisely, the solution $\ue$ of the
Allen-Cahn equation tends to a step function taking the value $+1$
on one side of an moving interface, and $-1$ on the other side.
This sharp interface, which we will denote by $\Gamma_t$, obeys
the law of motion $V_n=-\kappa$, where $V_n$ is the normal
velocity of $\Gamma _t$ in the exterior direction and $\kappa$ the
mean curvature at each point of $\Gamma_t$.

Then, some rigorous justification of this procedure were obtained.
In the framework of classical solutions, let us mention the works
of Bronsard and Kohn \cite{BS}, X. Chen \cite{C1,C2}, and de
Mottoni and Schatzman \cite{MS1, MS2}. Later, in \cite{AHM}, the
authors prove an optimal estimate for this convergence for
solutions with general initial data. By performing an analysis of
both the generation and the motion of interface, they show that
the solution develops a steep transition layer within a very short
time, and that the layer obeys the law of motion that coincides
with the formal asymptotic limit $V_n=-\kappa$ within an error
margin of order $\ep$ (previously, the best thickness estimate in
the literature was of order $\ep |\ln \ep|$, \cite{C1}). For
similar estimates of the thickness of the interface in related
problems we refer to \cite{A} (reaction-diffusion-convection
system as a model for chemotaxis with growth), \cite{HKMN}
(inhomogeneous Lotka-Volterra competition-diffusion system),
\cite{AGHMS} (fully anisotropic Allen-Cahn equation).

Since the classical motion by mean curvature may develop
singularities in finite time (extinction, \lq\lq pinch off"
phenomena...), one has to define a generalized motion by mean
curvature in order to study the singular limit of the Allen-Cahn
equation for all time. One represents $\Gamma _t$ as the level set
of an auxiliary function which solves (in the viscosity sense) a
nonlinear partial differential equation. This direct partial
differential equation approach was developed by Evans and Spruck
\cite{ES1}, Chen, Giga and Goto \cite{CGG}. In this framework of
viscosity solutions, we refer to Evans, Soner and Souganidis
\cite{ESS}, Barles, Soner and Souganidis \cite{BSS}, Barles ans
Souganidis \cite{BaS}, Ilmanen \cite{Ilm} for the singular limit
of reaction-diffusion equations, for all time.

\vskip 8 pt We now turn back to the nonlocal Allen-Cahn equation.
Problem $\Pe$ was considered by Chen, Hilhorst and Logak
\cite{CHL}. In order to underline its relevance in population
genetics and nervous transmission, they first show that $\Pe$ can
be seen as the limit, as $\sigma \to 0$ and $\tau \to 0$, of the
FitzHugh-Nagumo system
\[
 \quad\begin{cases}
 u_t&=\Delta u+\edeux f(u,\ep \frac{|\Omega|}\gamma v) \\
 \tau v_t&=\di \frac 1 \sigma \Delta v+u-\frac 1 \gamma v.
 \end{cases}
\]
Then, they study the motion of transition layers for the solutions
$\ue$ of $\Pe$. More precisely, for \lq\lq well-prepared" initial
data, they prove that, as $\ep \to 0$, the sharp interface limit,
which we will denote by $\Gamma_t$, obeys the law of motion
\[
 (P^0)\quad\begin{cases}
 \, V_{n}=-\kappa + c_0 (|\om _t ^+|-|\om _t ^-|)
 \quad \text { on } \Gamma_t \vspace{3pt}\\
 \, \Gamma_t\big|_{t=0}=\Gamma_0,
\end{cases}
\]
where $V_n$ is the normal velocity of $\Gamma _t$ in the exterior
direction, $\kappa$ the mean curvature at each point of
$\Gamma_t$,  $\om^-_t$ the region enclosed by $\Gamma_t$,
$\om^+_t$ the region enclosed between $\partial \om$ and
$\Gamma_t$, $c_0$ the constant defined by
\begin{equation}\label{intrinsic}
c_0=- \frac{\displaystyle \int _{-1}^{+1} \frac{\partial
f}{\partial v}(u,0)\,du}{\displaystyle \int _{-1}^{+1}
[2(W(u)-W(-1))]^{1/2}\,du},
\end{equation}
and $|A|$ the measure of the set $A$. As explained in \cite{CHL},
the Problem $(P^0)$ possesses a unique smooth solution locally in
time, say on some $[0,T]$. Moreover, in contrast with the \lq\lq
usual" motion by mean curvature which shrinks in finite time, the
nonlocal effect allows the possibility of nontrivial stationary
state (see \cite{Log} for a discussion in the radially symmetric
case).

\vskip 8 pt The goal of the present paper is to make a detailed
study of the limiting behavior of the solution $\ue$ of Problem
$\Pe$, without assuming that the initial datum already has a a
specific profile with a well-developed transition layer. In other
words, we study the generation of interface from arbitrary initial
data. Moreover, we obtain an improved error estimate, of $O(\ep)$,
between the solutions of $(P^\ep)$ and those of $(P^0)$.

Our main result, Theorem \ref{width}, describes the profile of the
solution after a very short initial period. It asserts that: given
a virtually arbitrary initial data $u_0$, the solution $\ue$
quickly becomes close to $\pm 1$, except in a small neighborhood
of the initial interface $\Gamma _0$, creating a steep transition
layer around $\Gamma _0$ ({\it generation of interface}). The time
needed to develop such a transition layer, which we will denote by
$t ^\ep$, is of order $\ep^2|\ln\ep|$. The theorem then states
that the solution $\ue$ remains close to the step function $\tilde
u$ on the time interval $[t^\ep,T]$ ({\it motion of interface}),
where $\tilde u$ is defined by
\begin{equation}\label{u}
\tilde u(x,t)=\begin{cases}
\, -1 &\text{ in } \om^-_t\\
\, +1 &\text{ in } \om^+_t
\end{cases} \quad\text{for } t\in[0,T].
\end{equation}
In other words, the motion of the transition layer is well
approximated by the limit interface equation $(P^0)$.

\begin{thm}[Generation, motion and thickness of transition layers]\label{width}
Let $\eta$ be an arbitrary constant satisfying $0< \eta <
\min(a+1,1-a)$ and set
$$
\mu=\tilde f'(a).
$$
Then there exist positive constants $\ep _0 $ and $\epaisseur$
such that, for all $\,\ep \in (0,\ep _0)$ and for all $\,t^\ep
\leq t \leq T$, where $t^\ep:=\mu ^{-1} \ep ^2 |\ln \ep|$, we have
\begin{equation}\label{resultat}
\ue(x,t) \in
\begin{cases}
\,[-1-\eta,+1+\eta]\quad\text{if}\quad
x\in\mathcal N_{\epaisseur \ep}(\Gamma_t)\\
\,[-1-\eta,-1+\eta]\quad\text{if}\quad
x\in\om_t^-\setminus\mathcal N_{\epaisseur\ep}(\Gamma
_t)\\
\,[+1-\eta,+1+\eta]\quad\text{if}\quad x\in\om
_t^+\setminus\mathcal N_{\epaisseur\ep}(\Gamma _t),
\end{cases}
\end{equation}
where $\mathcal N _r(\Gamma _t):=\{x\in \om, dist(x,\Gamma
_t)<r\}$ denotes the $r$-neighborhood of $\Gamma _t$.
\end{thm}

The estimate \eqref{resultat} implies that, once a transition
layer is formed, its thickness remains within order $\ep$ for the
rest of time.

\begin{cor}[Convergence]\label{total}
As $\ep\to 0$, $\ue$ converges to $\tilde u$ everywhere in $\cup
_{0<t\leq T}(\om^\pm_t\times\{ t\})$.
\end{cor}

This paper is organized as follows. In Section \ref{s:generation}
we study the generation of transition layers that takes place in a
very short time range. Section \ref{s:motion} is devoted to the
construction of a pair of sub- and super-solutions for the study
of the motion of interface. In Section \ref{s:proof} by fitting
the pair of sub- and super-solutions of Section \ref{s:generation}
into the pair of Section \ref{s:motion}, we prove our main result,
Theorem \ref{width}.  Since our arguments rely on a nonlocal
comparison principle borrowed from \cite{CHL}, we recall it in a
short appendix.

\section{Generation of interface}\label{s:generation}

In this section, we investigate the generation of interface,
namely the rapid formation of internal layers that takes place in
a neighborhood of $\Gamma_0=\{x\in \om|\ u_0(x)=a\}$ within the
time span of order $\ep^2 |\ln\ep|$. In this earlier stage, the
diffusion term is negligible and the partial differential equation
is approximated by the nonlocal equation $u_t = \edeux f(u,\ep
\int _\Omega u)$ and so, by the ordinary differential equation
\begin{equation}\label{ode-appro}
u_t = \edeux(\tilde
f(u)+O(\ep)).
\end{equation}

In the sequel, $\eta _0$ will stand for the quantity
$$
\eta _0:= \min (a+1,1 -a).
$$
The main result of the present section is the following.

\begin{thm}[Generation of interface]\label{th-gen}
Let $\eta \in (0,\eta _0)$ be arbitrary and define $\mu$ as the
derivative of $\tilde f(u)$ at the unstable zero $u=a$, that is
\begin{equation}\label{g-def-mu}
\mu={\tilde f}'(a).
\end{equation}
Then there exist positive constants $\ep_0$ and $M_0$ such that,
for all $\,\ep \in (0,\ep _0)$,
\begin{enumerate}
\item for all $x\in\om$,
\begin{equation}\label{g-part1}
-1-\eta \leq u^\ep(x,\mu ^{-1} \ep ^2 | \ln \ep |) \leq 1+\eta,
\end{equation}
\item for all $x\in\om$ such that $|u_0(x)-a|\geq M_0 \ep$, we
have that
\begin{align}
&\text{if}\;~~u_0(x)\geq a+M_0\ep\;~~\text{then}\;~~u^\ep(x,\mu
^{-1} \ep ^2 | \ln \ep |)
\geq 1-\eta,\label{g-part2}\\
&\text{if}\;~~u_0(x)\leq a-M_0\ep\;~~\text{then}\;~~u^\ep(x,\mu
^{-1} \ep ^2 | \ln \ep |)\leq -1+\eta \label{g-part3}.
\end{align}
\end{enumerate}
\end{thm}

The above theorem will be proved by constructing a suitable pair
of sub- and super-solutions based upon the ordinary differential
equation \eqref{ode-appro}. Note that the assumption of balanced
nonlinearity \eqref{int-f} is useless for the proof of the
generation of interface property.

\subsection{The bistable ordinary differential equation}

We first consider a slightly perturbed nonlinearity:
$$
\f(u):=\tilde f(u)+\delta ,
$$
where $\de$ is any constant. For $|\delta|$ small enough, this
function is still of the bistable type. More precisely, if $\de
_0$ is small enough, then for any $\de\in(-\de_0,\de_0)$, $\f$ has
exactly three zeros, namely $\hm < \h < \hp$, and there exists a
positive constant $C$ such that
\begin{equation}\label{h}
|\hm +1|+|\h -a|+|\hp-1|\leq C|\de|,
\end{equation}
\begin{equation}\label{mu}
| \mm -\mu| \leq C|\de|,
\end{equation}
where
\[
\mm:= {\f} '(\h)={\tilde f}'(\h).
\]

Now for each $\de\in(-\de_0,\de_0)$, we define $Y(\tau,\xi;\de)$
as the solution of the ordinary differential equation
\begin{equation}\label{ode}
\left\{\begin{array}{ll} Y_\tau (\tau,\xi;\de)&=\f
(Y(\tau,\xi;\de)) \quad \text { for }\tau >0\vspace{3pt}\\
Y(0,\xi;\de)&=\xi ,
\end{array}\right.
\end{equation}
where $\xi$ varies in $(-2C_0,2C_0)$, with $C_0$ being the
constant defined in \eqref{int1}. We claim that $Y(\tau,\xi;\de)$
has the following properties.

\begin{lem}\label{properties-Y} There exist positive
constants $\de_0$ and $C$ such that, for all $(\tau,\xi,\de) \in
(0,\infty)\times [-2C_0,2C_0]\times [-\de _0,\de _0]$,
\begin{enumerate}
 \item $|Y(\tau,\xi;\de)|\leq 2 C_0$
 \item $0<Y _ \xi (\tau,\xi;\de)$
  \item
$|\di{\frac{Y_{\xi\xi}}{Y_\xi}(\tau,\xi;\de)}|\leq C (e^{\mm
\tau}-1)$.
\end{enumerate}
\end{lem}
Property (i) is a direct consequence of the profile of $\f$ and so
of the qualitative properties of the solution of the bistable
ordinary differential equation \eqref{ode}; for proofs of (ii) and
(iii) we refer to \cite{AHM}, subsection 4.1.\qed

\subsection{Construction of sub- and super-solutions}

We are now ready to construct a pair of sub- and super-solutions in order to
prove the generation of interface property. By using some cut-off
initial data (see \cite{AHM}, subsection 3.2) we can modify
slightly $u_0$ near the boundary $\partial\Omega$ and make,
without loss of generality, the additional assumption
\begin{equation}\label{int2}
\frac{\partial u_0}{\partial \nu} =0 \quad\ \text{on}\; \
\partial\Omega.
\end{equation}
We set
\[
w_\ep^\pm(x,t)=Y\Big(\frac{t}{\ep^2},u_0(x)\pm\ep^2r(\pm \ep
\mathcal G,\frac{t}{\ep^2});\pm \ep \mathcal G\Big),
\]
where the function $r(\de,\tau)$ is given by
\[
r(\de,\tau)=\cste (e^{\mm \tau}-1),
\]
and the constant $\mathcal G$ by
$$
\mathcal G= 2C_0|\Omega| \max _{(u,v)\in[-2C_0,2C_0]\times[-1,1]}|
\frac{\partial f}{\partial v}(u,v)|.
$$

\begin{lem}\label{lemma-motion}
There exist positive constants $\ep_0$ and $\cste$ such that, for
all $\, \ep \in (0,\ep _0)$, $(w_\ep^-,w_\ep^+)$ is a pair of sub-
and super-solutions for Problem $\Pe$, in the domain $\om \times
(0,\mu ^{-1} \ep^2|\ln \ep|)$.
\end{lem}

Before proving the lemma, we remark that $w^-_\ep(x,0)=w^+ _\ep
(x,0)=u_0(x)$. Consequently, by the comparison principle, we obtain
\begin{equation}\label{coincee1}
w_\ep^-(x,t) \leq u^\ep(x,t) \leq w_\ep^+(x,t) \quad\ \text{ for
all } \ \ombar\times [0,\mu ^{-1} \ep^2|\ln \ep|].
\end{equation}

\noindent {\bf Proof.} First, the inequality $w_\ep
^-\leq w_\ep ^+$ follows from the fact that $Y(\tau,\xi;\de)$ increases
with both $\xi$ (see (ii) Lemma \ref{properties-Y}) and $\de$ (as easily seen from
the ordinary differential equation). Next, \eqref{int2} implies that both
$w_\ep ^+$ and $w_\ep ^-$ satisfy the Neumann homogeneous boundary
conditions. Hence, it remains to prove the inequalities $\mathcal
L _+ w_\ep ^+\geq 0$ and $\mathcal L _- w_\ep ^-\leq 0$ (see
Definition \ref{def-sub-sup}), provided that the constants $\ep
_0$ and $\cste$ are appropriately chosen.

If $\ep _0$ is sufficiently small, we note that $\pm \ep \mathcal
G \in (-\de _0,\de _0)$ and that, in the range $0 \leq t \leq \mu
^{-1} \ep ^2|\ln \ep|$,
$$
|\ep^2 \cste (e^{\mu(\pm \ep \mathcal G)t/\ep^2}-1)| \leq \ep
^2\cste (\ep^{-\mu(\pm \ep \mathcal G)/\mu}-1) \leq C_0,
$$
using \eqref{mu}. The above inequality implies
$$
u_0(x)\pm \ep^2r(\pm \ep \mathcal G,\frac{t}{\ep^2}) \in
[-2C_0,2C_0].
$$
These observations allow us to use the results of Lemma
\ref{properties-Y} with the choices $\tau:=t/\ep^2$,
$\xi:=u_0(x)\pm \ep^2r(\pm \ep \mathcal G,t/\ep ^2)$ and $\de:=\pm
\ep \mathcal G$. In particular, it follows from property (i) that
$ | \di \int _\Omega w^{\pm}_\ep (x,t)\, dx | \leq 2C_0 |\Omega|$
which in turn implies (thanks to the choice of $\mathcal G$)
\begin{equation}
\max _{\int _\Omega w_\ep ^- \leq s \leq \int _\Omega w_\ep ^+}
f(w_\ep ^+,\ep s)\leq \tilde f(w^+_\ep)+\ep \mathcal G.
\end{equation}

In view of the above inequality, some straightforward calculations
yield
$$
\mathcal L _+ w_\ep^+ \geq \edeux Y_\tau+\cste  \mu(\ep\mathcal
G)\emutt Y_\xi-|\nabla u_0|^2 Y_{\xi\xi}-\Delta u_0 Y_\xi -\edeux
\tilde f(Y)-\frac 1 \ep \mathcal G,
$$
where the argument $\left(\di \frac t {\ep^2},u_0(x)+\ep ^2 \cste (\emutt -1);\ep\mathcal G\right)$
of the function $Y$ and its derivatives is omitted. Noticing that the ordinary differential equation \eqref{ode}
writes as $Y_\tau =\tilde f(Y)+\ep\mathcal G$, we get
$$
\mathcal L _+ w_\ep^+ \geq
 Y_\xi\left[\cste  \mu(\ep
\mathcal G)\emutt-\Delta
u_0-\displaystyle{\frac{Y_{\xi\xi}}{Y_\xi}}|\nabla u_0|^2\right ].
$$
Using the estimate (iii) in Lemma \ref{properties-Y}, we obtain
\[
\begin{array}{ll}\mathcal L _+ w_\ep^+
&\geq Y_\xi\Big[\cste \mu(\ep \mathcal G)\emutt-|\Delta
u_0|-C(\emutt -1)|\nabla u_0|^2\Big] \vsp \\
&\geq Y_\xi\Big[(\cste \mu(\ep \mathcal G)-C |\nabla
u_0|^2)\emutt-|\Delta u_0| +C|\nabla
u_0|^2\Big].\\
\end{array}
\]
In view of \eqref{mu}, this inequality implies that, for $\ep \in
(0,\ep _0)$, with $\ep _0$ small enough,
\[
\mathcal L _+ w_\ep^+ \geq Y_\xi \Big[\cste  \frac
12\mu-C{C_0}^2-C_0\Big] \geq 0,
\]
by choosing $\cste$ large enough.

 Since one can prove $\mathcal L _- w_\ep ^- \leq 0$ by
similar arguments, this completes the proof of Lemma
\ref{lemma-motion}.\qed

\subsection{Proof of the generation of interface property}

In order to prove Theorem \ref{th-gen} we first quote a lemma from
\cite{AHM}; it makes more precise the bistable behavior of the ordinary differential equation
by giving  basic estimates of the function
$Y(\tau,\xi;\pm \ep \mathcal G)$ at time $\tau= \mu ^{-1}
|\ln \ep|$.

\begin{lem}\label{after-time}
Let $\eta \in (0,\eta _0)$ be arbitrary; there exist positive
constants $\ep_0$ and $\autrecste$ such that, for all
$\,\ep\in(0,\ep _0)$,
\begin{enumerate}
\item for all $\xi\in (-2C_0,2C_0)$,
\begin{equation}\label{part11}
-1-\eta \leq Y(\mu ^{-1} | \ln \ep |,\xi;\pm \ep \mathcal G) \leq
1+\eta;
\end{equation}
\item for all $\xi\in (-2C_0,2C_0)$ such that $|\xi-a|\geq
\autrecste \ep$, we have that
\begin{align}
&\text{if}\;~~\xi\geq a+\autrecste \ep\;~~\text{then}\;~~Y(\mu
^{-1}| \ln \ep |,\xi;\pm \ep \mathcal G)
\geq 1-\eta,\label{part22}\\
&\text{if}\;~~\xi\leq a-\autrecste \ep\;~~\text{then}\;~~Y(\mu
^{-1}| \ln \ep |,\xi;\pm \ep \mathcal G)\leq -1+\eta
\label{part33}.
\end{align}
\end{enumerate}
\end{lem}

\vskip 8pt \noindent{\bf Proof of Theorem \ref{th-gen}.} By
setting $t=\mu ^{-1} \ep ^2|\ln \ep|$ in \eqref{coincee1}, we get
\begin{multline}\label{gr}
Y\Big(\mu ^{-1}|\ln \ep|, u_0(x)-\ep^2 r(-\ep \mathcal G, \mu
^{-1}|\ln \ep|);-\ep \mathcal
G\Big)\\
\leq u^\ep(x,\mu ^{-1} \ep^2|\ln \ep|) \leq Y\Big(\mu ^{-1}|\ln
\ep|, u_0(x)+\ep^2 r(\ep \mathcal G, \mu ^{-1}|\ln \ep|);+\ep
\mathcal G \Big).
\end{multline}
We note that \eqref{mu} implies
\begin{equation}\label{point}
\lim _{\ep \rightarrow 0} \frac{\mu-\mu(\pm\ep \mathcal
G)}{\mu}\ln \ep=0,
\end{equation}
so that, if $\ep _0$ is sufficiently small,
$$
\ep ^2r(\pm \ep \mathcal G,\mu ^{-1} |\ln \ep|)= \cste\ep(\ep
^{(\mu-\mu(\pm\ep \mathcal G))/\mu}-\ep) \in (\frac 1 2 \cste\ep,
\frac 3 2 \cste \ep),
$$
and, for all $x \in \Omega$, it holds that $u_0(x)\pm \ep ^2r(\pm
\ep \mathcal G,\mu ^{-1} |\ln \ep|) \in (-2C_0,2C_0)$. Hence, the
result \eqref{g-part1} of Theorem \ref{th-gen} is a direct
consequence of \eqref{part11} and \eqref{gr}.

Next we prove \eqref{g-part2}. We take $x\in \om$ such that
$u_0(x)\geq a+M_0 \ep$; then
$$
\begin{array}{ll}u_0(x)-\ep^2r(-\ep \mathcal G,
\mu ^{-1}|\ln \ep|)
&\geq a+M_0\ep-\frac 3 2 \cste \ep\vspace{3pt}\\
&\geq a+\autrecste \ep,
\end{array}
$$
if we choose $M_0$ large enough. Using \eqref{gr} and
\eqref{part22} we see that inequality \eqref{g-part2} is true. The
inequality \eqref{g-part3} can be shown the same way. This
completes the proof of Theorem \ref{th-gen}.\qed

\section{Motion of interface}\label{s:motion}

In Section \ref{s:generation}, we have proved that the solution
$u^\ep$ of Problem $(P ^\ep)$ develops a clear transition layer
within a very short time.  The aim of the present section is to
show that, once such a clear transition layer is formed, it
persists for the rest of time and that its law of motion is well
approximated by the interface equation $(P^0)$. In order to study
this latter time range where the motion of interface occurs, we
will construct another pair of sub- and super-solutions
$(u_\ep^-,u_\ep^+)$ for Problem $\Pe$. To begin with we present
mathematical tools which are essential for this construction.

\subsection{Preliminaries}

{\bf The \lq \lq cut-off signed distance function\rq\rq.} Let
$\Gamma=\cup _{0< t \leq T}(\Gamma_t\times\{t\})$ be the solution
of the limit geometric motion problem $(P^0)$ and let $\widetilde d$ be
the signed distance function to $\Gamma$ defined by:
\begin{equation}\label{eq:dist}
\widetilde d (x,t)=
\begin{cases}
-&\hspace{-10pt}\mbox{dist}(x,\Gamma _t)\quad\text{for }x\in\Omega _t ^- \vspace{3pt}\\
&\hspace{-10pt} \mbox{dist}(x,\Gamma _t) \quad \text{for }
x\in\Omega _t ^+ ,
\end{cases}
\end{equation}
where $\mbox{dist}(x,\Gamma _t)$ is the distance from $x$ to the
hypersurface  $\Gamma _t$ in $\om$. The \lq\lq cut-off
signed distance function" $d$ is defined as follows. First,
choose $d_0>0$ small enough so that the signed distance function
$\widetilde d$ defined in \eqref{eq:dist} is smooth in the
following tubular neighborhood of $\Gamma$:
\[
 \{(x,t) \in \bar \Omega \times [0,T]\ |\ |\widetilde{d}(x,t)|<3d_0\},
\]
and that
\begin{equation}\label{front}
 \mbox{dist}(\Gamma_t,\partial \Omega)\geq 3d_0 \quad \textrm{ for all }
 t\in[0,T].
\end{equation}
Next let $\zeta(s)$ be a smooth increasing function on $\R$ such
that
\[
 \zeta(s)= \left\{\begin{array}{ll}
 s &\textrm{ if }\ |s| \leq d_0\vspace{4pt}\\
 -2d_0 &\textrm{ if } \ s \leq -2d_0\vspace{4pt}\\
 2d_0 &\textrm{ if } \ s \geq 2d_0.
 \end{array}\right.
\]
We then define the cut-off signed distance function $d$ by
\begin{equation}
d(x,t)=\zeta\left(\tilde{d}(x,t)\right).
\end{equation}
Note that, in view of
\eqref{front} and the definition of $d$, the equality $\nabla d=0$ holds in a
neighborhood of $\partial \Omega$. Note also that the equality $|\nabla d|=1$ holds in the region $\{(x,t) \in
\bar \Omega \times [0,T]\ |\ |d(x,t)|<d_0\}$. Moreover, since $\n d$ coincides with the outward normal unit vector to the hypersurface $\Gamma _t$, we
have $d_{t}(x,t)=-V_n$, where $V_n$ is the normal
velocity of the interface $\Gamma _t$ in the exterior direction. It is also known that the
mean curvature $\kappa$ of the interface is equal to $\Delta
 d$. Hence, since the moving interface
$\Gamma$ satisfies Problem $(P^0)$, an alternative equation for
$\Gamma$ is given by
\begin{equation}\label{FBP-aniso}
d _t=\Delta d -c_0 \nonlocal (t) \quad \text{ on }\Gamma _t,
\end{equation}
where $\gamma (t):=|\om _t ^+|-|\om _t ^-|$.

\vskip 8 pt \noindent {\bf The one dimensional standing wave
$U_0$.} Let $U_0(z)$ be the unique solution of the stationary
problem
\begin{equation}\label{eq-phi}
\left\{\begin{array}{ll}
{U_0} '' +\tilde f(U_0)=0 \vspace{3pt}\\
U_0(-\infty)= -1,\quad U_0(0)=a,\quad U_0(+\infty)=+1 .
\end{array} \right.
\end{equation}
This solution represents the first approximation of the profile of
a transition layer around the interface observed in the stretched
coordinates; it naturally arises when performing a formal
asymptotic expansion of the solution $\ue$ (see \cite{AHM},
Section 2). Note that the \lq\lq balanced stability
assumption\rq\rq, i.e. the integral condition $\di \int
_{-1}^{+1}\tilde f(u)\,du =0$, guarantees the existence of such a
standing wave. In the simple case where $\tilde f(u)=u(1-u^2)$, we
know that $U_0(z)=\tanh (z /\sqrt 2)$. In the general case, the
following standard estimates hold.
\begin{lem}\label{est-phi}
There exist positive constants $C$ and $\lambda$ such that
$$
\begin{array}{ll}
0 <1-U_0(z)&\leq Ce^{-\lambda|z|} \quad \text{ for } z\geq 0
\vspace{3pt}\\
0 <U_0(z)+1&\leq Ce^{-\lambda|z|} \quad \text{ for } z\leq 0.\\
\end{array}
$$
In addition, $U_0$ is a strictly increasing function and, for
$j=1, 2$,
\begin{equation}\label{diff-U0}
|D^jU_0(z)|\leq Ce^{-\lambda|z|} \quad \text{ for } z\in \R.
\end{equation}
\end{lem}

\vskip 8pt \noindent {\bf The solution $U_1$ of a linearized
problem.} Let $U_1(z,t)$ be the solution of the problem
\begin{equation}\label{eqU1}
\left\{\begin{array}{ll}
U_{1zz}+{\tilde f}'(U_0(z))U_1=\nonlocal(t)\left( -\displaystyle \frac{\partial f}{\partial v}(U_0(z),0)-c_0 {U_0}'(z)\right),\vsp\\
U_1(0,t)=0, \qquad\quad U_1(\cdot,t) \in L^\infty(\R),
\end{array}\right.
\end{equation}
where
\begin{equation}\label{cste-nonlocal}
c_0:=- \frac{\displaystyle \int _{\R} \frac{\partial f}{\partial
v}(U_0(z),0) {U_0}'(z)\,dz}{\displaystyle \int _{\R}
\left({U_0}'\right)^2(z)\,dz}.
\end{equation}
Again, the above problem
arises when performing a formal asymptotic expansion of the
solution $\ue$. Since \eqref{eqU1} can be seen as a linearized
problem for \eqref{eq-phi}, its solvability follows from a
Fredholm alternative: thanks to the definition of $c_0$, ${U_0}'$
turns out to be orthogonal to the right-hand side member of
\eqref{eqU1}. Moreover, there exist constants $M>0$ and $C>0$ such
that
\begin{equation}\label{def-M}
|U_1(z,t)|\leq M,
\end{equation}
\begin{equation}\label{delta-U1}
|U_{1t}(z,t)|\leq C,
\end{equation}
\begin{equation}\label{est-psi}
|U_{1z}(z,t)|+|U_{1zz}(z,t)| \leq Ce^{-\lambda |z|},
\end{equation}
for all $(z,t)\in\R\times[0,T]$. We omit the details and refer the
reader to \cite{AHM}, Section 2.

Note that, by multiplying equation \eqref{eq-phi} by ${U_0}'$ and integrating from
$-\infty$ to $z$, we obtain ${U_0}'(z)=[2(W(U_0(z))-W(-1))]^{1/2}$. Using this, it is now a matter
of routine to deduce from \eqref{cste-nonlocal} the more intrinsic expression \eqref{intrinsic}.

\subsection{Construction of sub- and super-solutions}

We look for a pair of sub- and super-solutions $u_\ep^{\pm}$ for
$(P^\ep)$ of the form
\begin{equation}\label{sub}
u_\ep^{\pm}(x,t)=U_0\left(\frac{d(x,t) \pm \ep p(t)}{\ep}\right)+\ep
U_1 \left(\frac{d(x,t) \pm \ep p(t)}{\ep},t\right)\pm q(t),
\end{equation}
where
\[
\begin{array}{lll}
p(t)=-\EB+e^{Lt}+ K, \\
q(t)=\sigma \big( \beta \EB+\ep^2 Le^{Lt}\big).
\end{array}
\]
Note that $q=\sigma\ep^2\,p_t$.  It is clear from the definition
of $u_\ep^\pm$ that
\begin{equation}\label{sub-lim}
\lim_{\ep\rightarrow 0} u_\ep^\pm(x,t)= \left\{
\begin{array}{ll}
+1 &\textrm { for all } (x,t) \in \cup_{\,0\leq t\leq T}(\Omega^+_t \times\{t\}) \vspace{4pt}\\
-1 &\textrm { for all } (x,t) \in \cup_{\,0\leq t\leq T}(\Omega^-_t \times\{t\}).\\
\end{array}\right.
\end{equation}

The main result of this section is the following:

\begin{lem}\label{fix}
Choose $\beta>0$ and $\sigma>0$ appropriately. Then for any $K>1$,
there exist positive constants $\ep_0$ and $L$ such that, for any
$\ep\in(0,\ep _0)$, $(u_\ep^-,u_\ep^+)$ is a pair of sub- and
super-solutions for $(P^\ep)$ in the domain $\bar \Omega \times
[0,T]$.
\end{lem}

{\noindent \bf Proof.} First, we claim that \eqref{temps-initial} and \eqref{ordre2} hold
as a consequence of \eqref{uep-H}. Then, since $d$ is constant in a neighborhood
of $\partial \Omega$, both $u_\ep^+$ and $u_\ep ^-$ satisfy the
Neumann homogeneous boundary conditions. Hence it remains to prove
the inequalities $\mathcal L _ + u_\ep ^+\geq 0$ and $\mathcal L_
-u_\ep ^-\leq 0$, provided that the various constants are
appropriately chosen.

\vskip 8 pt

We start with some useful inequalities.  On the one hand, by
assumption \eqref{der-f}, there exist positive constants $b,\,m$
such that
\begin{equation}\label{bords}
\tilde f'(U_0(z))\leq -m \qquad \hbox{if} \quad U_0(z)\in
[-1,\,-1+b]\cup[1-b,\,1].
\end{equation}
On the other hand, since the region $\{z\in\R\,|\,U_0(z)\in
[-1+b,\,1-b] \,\}$ is compact and since ${U_0}'>0$ on $\R$,
there exists a constant $a_1>0$ such that
\begin{equation}\label{milieu}
{U_0}'(z) \geq a_1 \qquad\hbox{if} \quad U_0(z)\in
[-1+b,\,1-b].
\end{equation}
We set
\begin{equation}\label{beta}
\beta= \frac{m}{4}\,,
\end{equation}
and choose $\sigma$ that satisfies
\begin{equation}\label{sigma}
0< \sigma \leq \min\,(\sigma_0,\sigma_1,\sigma_2),
\end{equation}
where
\[
\sigma_0:=\frac{a_1}{\di m+F_1},\quad
\sigma_1:=\frac{1}{\beta+1},\quad \sigma
_2:=\frac{4\beta}{H(\beta+1)},
\]
the constants $F_1$ and $H$ being given by
\begin{equation}\label{hess}
F_1:=\Vert \tilde f'\Vert _{L^\infty(-1,1)},\;\; H:=\max
_{(u,v)\in[-3,3]\times[-1,1]} \Vert Hess _{(u,v)} f\Vert,
\end{equation}
where $\Vert A\Vert :=\max _{i,j}|a_{ij}|$. Combining \eqref{bords} and \eqref{milieu}, and considering that
$\sigma \leq \sigma _0$, we obtain
\begin{equation}\label{U0-f}
{U_0}'(z)-\sigma \tilde f'(U_0(z))\geq \sigma m \qquad \hbox{for} \ \
-\infty<z<\infty.
\end{equation}

\vskip 8 pt

Now let $K>1$ be arbitrary. In what follows we will show that
$\mathcal L _+ u _\ep ^+ \geq 0$ provided that the constants
$\ep_0$ and $L$ are appropriately chosen. We recall that $-1
<U_0<1$ and that $|U_1|\leq M$. We go on under the following
assumption
\begin{equation}\label{ep0M}
\ep_0M\leq 1, \qquad \ep _0^2 Le^{LT} \leq 1\, .
\end{equation}
Then, given any $\ep\in(0,\ep_0)$, since $\sigma \leq \sigma _1$
we have $0\leq q(t)\leq 1$, so that
\begin{equation}\label{uep-pm}
-3\leq u_\ep^\pm(x,t) \leq +3\, .
\end{equation}

\vskip 8 pt

In order to evaluate the \lq\lq nonlocal part" of $\mathcal L _+ u
_\ep ^+$, we need bounds for the quantities $\di \int _\Omega
u^\pm _\ep (x,t)\,dx$. For the sake of clarity, the arguments of
most of the functions are omitted in the following. We write $\di
\int _\Omega u_\ep ^+\,dx =\int _{\Omega _t ^+}(U_0 -1)\,dx +
|\Omega _t ^+|+\int _{\Omega _t ^-}(U_0 +1)\,dx - |\Omega _t
^-|+\int _\Omega(\ep U_1+q)\,dx$, which yields
$$
\begin{array}{ll}
\di{ \int _ \Omega u_\ep ^+\, dx-\nonlocal(t)}&=\di{|\Omega|q (t)+\int _\Omega \ep U_1\,dx}\vsp\\
&\;\;\;\;\;\;\;\;\;\; \di{+\int _{\Omega _t ^+}(U_0 -1)\,dx+\int _{\Omega _t ^-}(U_0 +1)\,
dx}\vsp\\
&=:|\Omega|q (t)+I_\ep(t)+I_+(t)+I_-(t).
\end{array}
$$
In the following we will denote by $C$ various positive constants
that are independent of $\ep \in (0,\ep _0)$. Since $U_1$ is
bounded, we have $|I_\ep (t)|\leq C\ep$.

In order to estimate $I_+(t)$, we use the partition
$$
\Omega _t
^+=\{x|\; d(x,t)\geq d_0\} \cup \{x|\; 0<d(x,t)<d_0\}.
$$
First assume $d(x,t)\geq d_0$. From Lemma \ref{est-phi}, we deduce
that
$$
0\leq 1-U_0(\frac{d(x,t)+\ep p(t)}\ep) \leq C e^{-\lambda|d(x,t)+\ep p(t)|/\ep}\leq Ce^{-
\lambda d_0/\ep},
$$
from which we infer that
$$
0\leq \int _ {d(x,t)\geq d_0} \left(1-U_0(\frac{d(x,t)+\ep
p(t)}\ep)\right)\,dx\leq Ce^{-\lambda d_0/\ep}.
$$
In order to estimate the integral on $\{x|\; 0<d(x,t)<d_0\}$, we
use arguments similar to those used in \cite{CHL}. We denote by
$J(s,d)$ the Jacobi of the transformation $x\mapsto(s,d)$, where
$s(x,t)$ is the projection of $x$ on $\Gamma _t$ along the normal
of $\Gamma _t$ and $d(x,t)(=\tilde d (x,t))$ is the signed
distance defined above; we define $C_J:=\max_{0\leq t\leq T} \Vert
J(\cdot,\cdot)\Vert_{L ^\infty (\Gamma _t \times [-d_0,d_0])}$.
This yields
$$
\begin{array}{lll}
0\leq \di{\int _{0<d(x,t)< d_0}}&\left(1-U_0(\di{\frac{d(x,t)+\ep p(t)}\ep})\right)\,dx\vsp\\
&=\di{\int _{\Gamma _t}\int _0 ^{d_0}}\left(1-U_0(
\frac{r+\ep p(t)}\ep)\right)J(s,r)\,drds\vsp\\
&\leq C_J \ep\di{\int _{\Gamma _t}\int _0 ^{d_0/\ep}} \left(1-U_0(z+p(t))\right)\,dzds\vsp\\
&\leq C_J\ep\di{\int _{\Gamma _t}\int _0
^{+\infty}}(1-U_0(u))\,duds\leq C\ep .
\end{array}
$$

As far as $I_-(t)$ is concerned, we first assume that $d(x,t)\leq
-d_0$. Note that $0<K-1 \leq p \leq e^{LT} +K$. Consequently, if
we assume
\begin{equation}\label{ga}
e^{LT}+K \leq \frac{d_0}{2\ep_0},
\end{equation}
then $\displaystyle{\frac{d_0}{\ep}}-|p|\geq
\displaystyle{\frac{d_0}{2\ep}}$. By using similar arguments as
the ones above, we obtain
$$
0\leq \int _ {d(x,t)\leq -d_0} \left(U_0(\frac{d(x,t)+\ep
p(t)}\ep)+1\right)\,dx\leq Ce^{-\lambda d_0/(2\ep)}.
$$
Concerning the region $\{x|\; -d_0<d(x,t)<0\}$, we get
$$
\begin{array}{lll}
0\leq \di{\int _{-d_0<d(x,t)< 0}}&\left(U_0(\di{\frac{d(x,t)+\ep p(t)}\ep})+1\right)\,dx\vsp\\
&=\di{\int _{\Gamma _t}\int _{-d_0} ^0}\left(U_0(
\frac{r+\ep p(t)}\ep)+1\right)J(s,r)\,drds\vsp\\
&\leq C_J \ep\di{\int _{\Gamma _t}\int _{-d_0/\ep} ^0 } \left(U_0(z+p(t))+1)\right)\,dzds\vsp\\
&\leq C_J\ep\di{\int _{\Gamma _t}\int _{-\infty}
^{d_0/(2\ep_0)}}(U_0(u)+1)\,duds\leq C\ep.
\end{array}
$$

Since one would obtain similar estimates with $u_\ep ^+$ replaced by  $u_\ep ^-$, the above estimates yield
\begin{equation}
\left | \int _\Omega u^{\pm}_\ep\, dx-\nonlocal(t)\right | \leq
C\ep + Cq(t),
\end{equation}
which, in turn, implies
\begin{equation}
\max _{\int _\Omega u_\ep ^- \leq s \leq \int _\Omega u_\ep ^+}
f(u_\ep ^+,\ep s)\leq f(u^+_\ep,\ep \nonlocal(t))+ C \ep ^2 +C\ep
q(t)
\end{equation}
and
\begin{equation}
\min _{\int _\Omega u_\ep ^- \leq s \leq \int _\Omega u_\ep ^+}
f(u_\ep ^-,\ep s)\geq f(u^-_\ep,\ep \nonlocal(t))- C \ep ^2 -C\ep
q(t).
\end{equation}

\vskip 8 pt

Now, we can turn back to the proof of $\mathcal L _+ u_\ep ^+\geq 0$. From the above inequality, we get
\begin{equation}
\mathcal L _+ u^+_\ep \geq (u_\ep ^+)_t-\Delta u_\ep ^+ -\edeux
f(u^+_\ep,\ep\nonlocal (t))-C-C\frac 1 \ep q(t).
\end{equation}
Straightforward computations yield
$$
\begin{array}{lll}
(u_\ep^+)_t= {U_0}'(\displaystyle{\frac{d_t}{\ep}}+p_t) +\ep
U_{1t}+ U_{1z}(d_t+\ep p_t) +q_t\vsp \\
\Delta u_\ep^+= {U_0}''\displaystyle{\frac{|\nabla d|^2}{\ep ^2}}
+ {U_0}'\displaystyle{\frac{\Delta d}{\ep}} + U_{1zz}
\displaystyle{\frac{|\nabla d|^2}{\ep}} + U_{1z} \Delta  d,
\end{array}
$$
where the function $U_0$, as well as its derivatives, are
evaluated at $z=\big(d (x,t)+\ep p(t)\big)/ \ep $, whereas the
function $U_1$, as well as its derivatives, are evaluated at
$\Big(\big(d (x,t)+\ep p(t)\big)/ \ep,t \Big)$. We also have
\begin{align*}
f(u_\ep ^+,\ep \nonlocal(t))\leq &\tilde f(U_0)+(\ep U_1+q)\tilde
f '(U_0)+\ep \nonlocal(t) \frac{\partial f}{\partial
v}(U_0,0)\\
&+H\left(\frac 12(\ep U_1+q)^2+\frac 12 \ep ^2 \nonlocal ^2(t)+
(\ep U_1+q)\ep \nonlocal(t)\right),
\end{align*}
where $H$ was defined in \eqref{hess}. Combining the above
expressions with the equations \eqref{eq-phi} and \eqref{eqU1} for
$U_0$ and $U_1$, we obtain
\[
\mathcal L _+u_\ep^+\geq E_1+\cdots+E_6,
\]
where:
\vsp \\
$\qquad\quad  E_1=- \edeux q\,\left(\tilde f'(U_0)
+\frac 12  H q\right)+{U_0}'p_t+q_t$\vsp \\
$\qquad\quad  E_2=\displaystyle{ \left(\frac{{U_0}''}{\ep^2}
+ \frac {U_{1zz}}{\ep}\right)}(1-|\n d|^2)$\vsp \\
$\qquad\quad  E_3=\displaystyle{ \left(\frac {{U_0}'}{\ep}+
U_{1z}\right)}(d_t-\Delta d +c_0\nonlocal(t))$\vsp \\
$\qquad\quad  E_4=\ep U_{1z} p_t+\di{\frac{1}{\ep}}
q\,\left(-H\nonlocal(t)-H U_1-C \right)$\vsp\\
$\qquad\quad  E_5=-c_0 \nonlocal  (t)U_{1z}-\frac 12 H{U_1}^2
-\frac
12H \nonlocal ^2(t)-H \nonlocal (t)U_1-C$\vsp \\
$\qquad\quad  E_6=\ep U_{1t}$.\\

In the sequel, we estimate the terms $E_1$---$E_6$ and denote by
$C_i$ various positive constants that are independent of $\ep$.

\subsubsection {The term $E_1$}

Direct computation gives
$$
E_1=\frac{\beta}{\ep^2}\,\EB(I-\sigma\beta)+Le^{Lt}(I+\ep^2\sigma
L),
$$
where
$$
I={U_0}'-\sigma \tilde f '(U_0)-\frac {\sigma^2}2 H(\beta\EB+\ep^2
Le^{Lt}).
$$
In virtue of \eqref{U0-f}, we have
\[
I\geq \sigma m-\frac {\sigma^2}{2} H(\beta+\ep^2 Le^{LT}).
\]
Combining this, \eqref{beta}, \eqref{ep0M} and the inequality
$\sigma \leq \sigma _2$, we obtain $ I \geq 2\sigma\beta$.
Consequently, we have
$$
E_1\geq \frac{\sigma\beta^2}{\ep^2}\EB + 2\sigma\beta L e^{Lt}.
$$

\subsubsection {The term $E_2$}

First, in the region where $|d|<d_0$, we have $|\n d|=1$, hence
$E_2=0$. Next we consider the region where $|d|\geq d_0.$ We
deduce from Lemma \ref{est-phi} and from \eqref{est-psi} that
$$
|E_2|\leq C(\frac{1}{\ep^2}+\frac{1}{\ep})e^{-\lambda|d+\ep p|/
\ep}\leq \frac{2C}{\ep^2}e^{-\lambda(d_0 / \ep-|p|)}.
$$
In view of \eqref{ga} we have $0< p \leq \di{\frac{d_0}{2\ep}}$ so
that
$$
|E_2|\leq \frac{2C}{\ep^2}e^{-\lambda d_0 / (2\ep)} \leq C_2.
$$

\subsubsection {The term $E_3$}

Recall that
$$
(d_t-\Delta d)(x,t) +c_0 \nonlocal(t)=0 \qquad \textrm{on} \quad
\Gamma_t=\{x \in \om,\; d(x,t)=0\}.
$$
Since the interface $\Gamma_t$ is smooth, both $\Delta d$ and
$d_t$ are Lipschitz continuous near $\Gamma_t$.  It follows from
the mean value theorem applied on both sides of $\Gamma _t$ that
there exists a constant $N>0$ such that:
$$
|(d_t-\Delta d)(x,t) +c_0\nonlocal (t)|\leq N|d(x,t)| \quad
\textrm{ for all }  (x,t) \in \Omega \times (0,T).
$$
Applying  Lemma \ref{est-phi} and the estimate \eqref{est-psi}
we deduce that
$$
\begin{array}{lll}
|E_3|&\leq 2NC\displaystyle{\frac
{|d|}{\ep}}e^{-\lambda| d/ \ep +p|}\vsp \\
&\leq 2NC \max_{\xi \in \R }|\xi|e^{-
\lambda|\xi +p|}\vsp \\
&\leq 2NC\max (|p|,\di{\frac 1 \lambda}).
\end{array}
$$
Thus, recalling that $|p|\leq e^{Lt}+K$, we obtain
$$
|E_3|\leq C_3(e^{Lt}+K)+{C_3}', $$ where $C_3:=2NC$ and
${C_3}':=2NC/\lambda$.

\subsubsection {The terms $E_4$, $E_5$ and $E_6$}

Since $U_1$, $U_{1z}$, $U_{1t}$ and $\nonlocal$ are bounded, it is
a matter of routine to see that
\[
|E_4| \leq C_4\big ( \frac 1 \ep \beta \EB+\ep L e^{Lt}\big ),\;\; |E_5|\leq C_5,\;\; |E_6|\leq \ep C_6.
\]

\subsubsection {Completion of the proof}

Collecting all these estimates gives
\begin{equation}\label{september}
 \mathcal L _+ u_\ep^+\geq (\frac{\sigma \beta ^2}{\ep^2}-\frac{C_4\beta}{\ep})
 e^{-\beta t/\ep ^2}+ (2\sigma \beta L-C_3-\ep C_4 L)e^{Lt}-C_7,
\end{equation}
where $ C_7:=C_2+KC_3+{C_3}'+C_5+C_6$. Now we set
\[
 L:=\frac 1 T\ln \frac {d_0}{4\ep _0},
\]
which, for $\ep_ 0$ small enough, validates assumptions
\eqref{ep0M} and \eqref{ga}. For $\ep_0$ small enough, the first
term of the right-hand side of \eqref{september} is positive, and
\[
 \mathcal L _+ u_\ep^+\geq \big[\sigma \beta L-C_3]e^{Lt}-C_7 \geq  \frac 1 2 \sigma \beta L -C_7 \geq 0.
\]
The proof of \eqref{super-solution} is now complete, with the choice of
the constants $\beta, \sigma$ as in \eqref{beta},
\eqref{sigma}. Since one can prove \eqref{sub-solution} by similar arguments, this completes the proof of Lemma \ref{fix}. \qed

\section{Proof of Theorem \ref{width}}\label{s:proof}

Let $\eta \in (0,\eta _0)$ be arbitrary. Choose $\beta$ and
$\sigma$ that satisfy \eqref{beta}, \eqref{sigma} and
\begin{equation}\label{eta}
\sigma \beta \leq \frac \eta 3.
\end{equation}
By the generation of interface property, Theorem \ref{th-gen},
there exist positive constants $\ep_0$ and $M_0$ such that
\eqref{g-part1}, \eqref{g-part2} and \eqref{g-part3} hold with the
constant $\eta$ replaced by $ \sigma \beta /2$. Since $\n u_0 \neq
0$ everywhere on $\Gamma _0=\{x\in\om| \; u_0(x)=a\}$ and since
$\Gamma _0$ is a compact hypersurface, we can find a positive
constant $M_1$ such that
\begin{equation}\label{corres}
\begin{array}{ll}\text { if } \quad d_0 (x) \geq \ M_1 \ep
&\text { then } \quad u_0(x) \geq a +M_0 \ep\vspace{3pt}\\
\text { if } \quad d_0 (x) \leq -M_1 \ep & \text { then } \quad
u_0(x) \leq a -M_0 \ep.
\end{array}
\end{equation}
Here $d_0(x):= d(x,0)$ denotes the signed distance function
associated with the hypersurface $\Gamma_0$. Now we define
functions $H^+(x), H^-(x)$ by
\[
\begin{array}{l}
H^+(x)=\left\{
\begin{array}{ll}
+1+\sigma\beta/2\quad\ &\hbox{if}\ \ d_0(x)> -M_1\ep\\
-1+\sigma\beta/2\quad\ &\hbox{if}\ \ d_0(x)\leq  -M_1\ep,
\end{array}\right.
\vsp\\
H^-(x)=\left\{
\begin{array}{ll}
+1-\sigma\beta/2\quad\ &\hbox{if}\ \ d_0(x)\geq \;M_1\ep\\
-1-\sigma\beta/2\quad\ &\hbox{if}\ \ d_0(x)<  \;M_1\ep.
\end{array}\right.
\end{array}
\]
Then from \eqref{g-part1}, \eqref{g-part2}, \eqref{g-part3} (with
$\eta$ replaced by $ \sigma \beta /2$) and \eqref{corres}, we see
that
\begin{equation}\label{H-u}
H^-(x) \,\leq\, u^\ep(x,\mu^{-1} \ep^2|\ln \ep|) \,\leq\,
H^+(x)\qquad \hbox{for}\ \ x\in\Omega.
\end{equation}

Next we fix a sufficiently large constant $K>1$ such that
\begin{equation}\label{K}
U_0(-M_1+K) \geq 1-\frac {\sigma \beta}{3} \quad \text { and }
\quad U_0(M_1-K) \leq -1+\frac {\sigma \beta}{3}.
\end{equation}
For this $K$, we choose $\ep _0$ and $L$  as in Lemma \ref{fix}.
We claim that
\begin{equation}\label{uep-H}
u_\ep^-(x,0)\leq H^-(x),\quad\ H^+(x)\leq u_\ep^+(x,0) \qquad
\hbox{for} \ \ x\in\Omega,
\end{equation}
with $(u_\ep ^-,u_\ep ^+)$ the pair of sub- and super-solutions
defined in \eqref{sub} for the study of the motion of interface.
We shall only prove the former inequality, as the proof of the
latter is virtually the same. Then it amounts to showing that
\begin{equation}\label{c3}
u_\ep ^- (x,0)=U_0\big(\frac {d_0(x)}{\ep}-K\big)+\ep U_1\big(\frac{d_0(x)}{\ep}-K,0\big)-\sigma (\beta+\ep ^2 L) \;\leq\; H^-(x).
\end{equation}
Recall that $|U_1|\leq M$. Therefore, by choosing $\ep_0$ small
enough so that $\ep_0 M \leq \sigma\beta/6$, we see that
$$
u_\ep ^- (x,0)  \leq\;U_0\big(\frac
{d_0(x)}{\ep}-K\big)-\frac{5}{6}\sigma\beta.
$$
In the range where $d_0(x) < M_1 \ep$, the fact that $U_0$ is an
increasing function and the second inequality in \eqref{K} imply
\[
U_0\big(\frac {d_0(x)}{\ep}-K\big)-\frac{5}{6}\sigma\beta \;\leq\;
U_0(M_1-K)-\frac 56 \sigma \beta \;\leq\; -1-\frac {\sigma
\beta}{2}\;=\;H^-(x).
\]
On the other hand, in the range where $d_0(x) \geq M_1 \ep$, we
have
\[
U_0\big(\frac {d_0(x)}{\ep}-K\big)-\frac{5}{6}\sigma\beta \;\leq\;
1-\frac{5}{6}\sigma\beta \;\leq\;H^-(x).
\]
This proves \eqref{c3}, hence \eqref{uep-H} is established.

Combining \eqref{H-u} and \eqref{uep-H}, we obtain
$$
u_\ep^-(x,0)\leq u^\ep(x,\mu ^{-1} \ep ^2|\ln \ep|) \leq
u_\ep^+(x,0).
$$
Since $(u_\ep^-,u_\ep^+)$ is a pair of sub- and super-solutions
for Problem $\Pe$, the comparison principle yields
\begin{equation}\label{ok}
u_\ep^-(x,t) \leq u^\ep (x,t+t^\ep) \leq u_\ep^+(x,t) \quad \text
{ for } 0 \leq t \leq T-t^\ep,
\end{equation}
where $t^\ep=\mu ^{-1} \ep ^2|\ln \ep|$. Note that, in view of
\eqref{sub-lim}, this is enough to prove Corollary \ref{total}.
Now let $\epaisseur$ be a positive constant such that
\begin{equation}\label{C}
U_0(\epaisseur-e^{LT}-K) \geq 1-\frac \eta 2 \quad \text { and }
\quad U_0(-\epaisseur+e^{LT}+K) \leq -1+\frac \eta 2.
\end{equation}
One then easily checks, using \eqref{ok} and \eqref{eta}, that,
for $\ep _0$ small enough, for $0\leq t \leq T-t^\ep$, we have
\begin{equation}\label{correspon}
\begin{array}{ll}\text { if } \quad d(x,t) \geq \epaisseur \ep &
\text { then } \quad
u^\ep (x,t+t^\ep) \geq 1 -\eta\vspace{3pt}\\
\text { if } \quad d(x,t) \leq -\epaisseur \ep & \text { then }
\quad u^\ep (x,t+t^\ep) \leq -1 +\eta,
\end{array}
\end{equation}
and
$$
u^\ep (x,t+t^\ep) \in [-1-\eta,1+\eta],
$$
which completes the proof of Theorem \ref{width}.\qed

\setcounter{section}{1}
\renewcommand{\thesection}{\Alph{section}}
\section*{Appendix - Comparison principle}

The definition of sub- and super-solutions is the one proposed by
Chen, Hilhorst and Logak \cite{CHL}. It involves simultaneously a
super- and a sub-solution.

\begin{defi}\label{def-sub-sup}
Let $(u_\ep ^-,u_ \ep ^+)$ be a pair of smooth functions defined on $\bar \Omega \times [0,T]$ and satisfying
\begin{equation}\label{ordre}
u_\ep ^-\leq u_ \ep ^+ \quad \text{ in }\; \bar \Omega \times
[0,T],
\end{equation}
and
\begin{equation}\label{ordre-neumann}
\frac{\partial u_\ep ^-}{\partial \nu} \leq 0 \leq \frac{\partial
u_\ep ^+}{\partial \nu} \quad\ \text{ on } \; \ \partial \om
\times (0,T).
\end{equation}
We say that $(u_\ep ^-,u_ \ep ^+)$ is a pair of sub- and
super-solutions for Problem $(P^\ep)$ if
\begin{equation} \label{super-solution}
\mathcal L _+ u_\ep ^+:=(u_\ep ^+)_t-\Delta u_\ep ^+ -\edeux \max
_{\int _\Omega u_\ep ^- \leq s \leq \int _\Omega u_\ep ^+} f(u_\ep
^+,\ep s) \geq 0 \;\;\text{ in } \Omega\times (0,T),
\end{equation}
\begin{equation}\label{sub-solution}
\mathcal L _- u_\ep ^-:= (u_\ep ^-)_t-\Delta u_\ep ^- -\edeux \min
_{\int _\Omega u_\ep ^- \leq s \leq \int _\Omega u_\ep ^+} f(u_\ep
^-,\ep s) \leq 0 \;\;\text{ in } \Omega\times (0,T).
\end{equation}
\end{defi}

As proved in \cite{CHL}, the following comparison principle holds.

\begin{prop}\label{comparison}
Let a pair of sub- and super-solutions be given. Assume that, for
all $x \in \om$,
\begin{equation}\label{temps-initial}
u_\ep ^-(x,0) \leq u_0(x) \leq u_\ep ^+(x,0).
\end{equation}
Then, if we denote by $u^\ep$ the solution of Problem $\Pe$, the
function $\ue$ satisfies
\begin{equation}\label{compar}
u_\ep ^-(x,t) \leq \ue (x,t) \leq u_\ep ^+(x,t),
\end{equation}
for all $(x,t)\in \Omega \times [0,T]$.
\end{prop}

As easily seen from the proof in \cite{CHL}, one could replace the assumption \eqref{ordre} by
the assumption \eqref{temps-initial} together with the condition that
\begin{equation}\label{ordre2}
\int _\Omega u_\ep ^- (x,t) \,dx \leq \int _\Omega u_\ep ^+(x,t)\, dx,
\end{equation}
for all $t\in[0,t_0]$ with $t_0>0$. More precisely if \eqref{temps-initial}, \eqref{ordre2},
\eqref{ordre-neumann}, \eqref{super-solution} and \eqref{sub-solution} hold, then the conclusion \eqref{compar} follows.


\begin{thebibliography}{ABCD}

\bibitem{A} M.~Alfaro, {\it The singular limit of a chemotaxis-growth system with general initial data},
Adv. Differential Equations  {\bf 11}  (2006),  no. 11,
1227--1260.

\bibitem{AGHMS} M.~Alfaro, H. Garcke, D.~Hilhorst, H.~Matano and
R. Sch\"a\-tzle, {\it Motion by anisotropic mean curvature as
sharp interface limit of an inhomogeneous and anisotropic
Allen-Cahn equation}, submitted to Proc. Roy. Soc. Edinburgh Sect.
A.

\bibitem{AHM} M.~Alfaro, D.~Hilhorst and H.~Matano,
{\it The singular limit of the Allen-Cahn equation and the
FitzHugh-Nagumo system}, J.~Differential Equations {\bf 245}
(2008), 505--565.

\bibitem{AC} S.~Allen and J.~Cahn,
{\it A microscopic theory for antiphase boundary motion and its
application to antiphase domain coarsening}, Acta Metallica {\bf
27} (1979), 1084--1095.

\bibitem{BSS} G.~Barles, H.~M.~Soner and P.~E.~Souganidis,
{\it Front propagation and phase field theory}, SIAM J.~Control
Optim. {\bf 31} (1993), 439--469.

\bibitem{BaS} G.~Barles and P.~E.~Souganidis, {\it A new approach to
front propagation problems : theory and applications},
Arch.~Rat.~Mech.~Anal. {\bf 141} (1998), 237--296.

\bibitem{BS} L.~Bronsard and R.~V.~Kohn, {\it Motion by mean
curvature as the singular limit of Ginzburg--Landau dynamics}, J.
Differential Equations {\bf 90} (1991), 211--237.

\bibitem{C1} X.~Chen,
{\it Generation and propagation of interfaces for
reaction-diffusion equations}, J.~Differential Equations {\bf 96}
(1992), 116--141.

\bibitem{C2} X.~Chen, {\it Generation and propagation of
interfaces for reaction-diffusion systems},
Trans.~Amer.~Math.~Soc. {\bf 334} (1992), 877--913.

\bibitem{CHL} X.~Chen, D.~Hilhorst and E.~Logak, {\it Asymptotic
behavior of solutions of an Allen-Cahn equation with a nonlocal
term}, Nonlinear Anal. {\bf 28} (1997), no. 7, 1283--1298.

\bibitem{CGG} Y.~G.~Chen, Y.~Giga and S.~Goto, {\it Uniqueness and existence of viscosity solutions of
generalized mean curvature flow equations}, J.~Diff.~Geom. {\bf
33} (1991), 749--786.

\bibitem{ESS} L.~C.~Evans, H.~M.~Soner and P.~E.~Souganidis,
{\it Phase transitions and generalized motion by mean curvature},
Comm.~Pure Appl.~Math. {\bf 45} (1992), 1097--1123.

\bibitem{ES1} L.~C.~Evans and J. Spruck,
{\it Motion of level sets by mean curvature I}, J. Differential
Geometry {\bf 33} (1991), 635--681.

\bibitem{KO} K.~Kawasaki and T.~Ohta,
{\it Kinetic drumhead model of interface I}, Progress of
Theoretical Physics {\bf 67} (1982) 147--163.

\bibitem{HKMN} D.~Hilhorst, G.~Karali, H.~Matano and K.~Nakashima,
{\it Singular limit of a spatially inhomogeneous Lotka-Volterra
competition-diffusion system}, Comm. Partial Differential
Equations  {\bf 32}  (2007), 879--933.

\bibitem{Ilm} T.~Ilmanen, {\it Elliptic regularization and partial
regularity for motion by mean curvature}, Memoirs of the American
Mathematical Society, {\bf 108} (1994).

\bibitem{Log}E.~Logak, {\it Singular limit of reaction-diffusion systems and modified
motion by mean curvature}, Proc. Roy. Soc. Edinburgh A {\bf 132}
(2002), no. 4, 951--973.

\bibitem{MS1} P.~de Mottoni and M.~Schatzman, {\it Development of
interfaces in $\R^n$}, Proc.~Roy. Soc.~Edinburgh {\bf 116A}
(1990), 207--220.

\bibitem{MS2} P.~de Mottoni and M.~Schatzman,
{\it Geometrical evolution of developed interfaces},
Trans.~Amer.~Math.~Soc.~{\bf 347} (1995), 1533--1589.


\end{thebibliography}
\end{document}